\documentclass[review,onefignum,onetabnum]{siamart190516}

\usepackage{lipsum}
\usepackage{amsfonts}
\usepackage{graphicx}
\usepackage{epstopdf}
\usepackage{algorithmic}
\ifpdf
  \DeclareGraphicsExtensions{.eps,.pdf,.png,.jpg}
\else
  \DeclareGraphicsExtensions{.eps}
\fi

\usepackage{enumitem}
\setlist[enumerate]{leftmargin=.5in}
\setlist[itemize]{leftmargin=.5in}
\usepackage{todonotes} 
\usepackage{amsmath,amssymb}
\usepackage{placeins}
\usepackage{cite}
\usepackage{float}

\newsiamremark{remark}{Remark}
\newsiamremark{hypothesis}{Hypothesis}
\crefname{hypothesis}{Hypothesis}{Hypotheses}
\newsiamthm{claim}{Claim}

\newcommand{\eps}{\varepsilon}

\headers{Canonical models for torus canards in elliptic bursters}{Emre Baspinar, Daniele Avitabile, and Mathieu Desroches}

\title{Canonical models for torus canards in elliptic bursters\thanks{\funding{This work was supported by Human Brain Project (HBP) funded from the European Union’s Horizon 2020 Framework Programme for Research and Innovation under the Specific Grant Agreement No. 785907 (Human Brain Project SGA2).}}}

\author{tba}
\author{Emre Baspinar\thanks{MathNeuro Team, Inria Sophia Antipolis M{\'e}diterran{\'e}e, France  (\email{emre.baspinar@inria.fr}).}
\and Daniele Avitabile\thanks{Department of Mathematics, Vrije Universiteit Amsterdam, Netherlands (\email{d.avitabile@vu.nl}).}
\and Mathieu Desroches\thanks{MathNeuro Team, Inria Sophia Antipolis M{\'e}diterran{\'e}e, France  (\email{mathieu.desroches@inria.fr}).}}

\usepackage{amsopn}

\ifpdf
\hypersetup{
  pdftitle={Canonical models for torus canards in elliptic bursters},
  pdfauthor={}
}
\fi

\graphicspath{{figures/}}

\begin{document}

\maketitle

\begin{abstract}
We revisit elliptic bursting dynamics from the viewpoint of torus canard solutions. We show that at the transition to and from elliptic burstings, classical or \textit{mixed-type} torus canards can appear, the difference between the two being the fast subsystem bifurcation that they approach, saddle-node of cycles for the former and subcritical Hopf for the latter. We first showcase such dynamics in a Wilson-Cowan type elliptic bursting model, then we consider minimal models for elliptic bursters in view of finding  transitions to and from bursting solutions via both kinds of torus canards. We first consider the canonical model proposed in~\cite{izhikevich2000subcritical} and adapted to elliptic bursting in~\cite{ju2018bottom}, and we show that it does not produce mixed-type torus canards due to a nongeneric transition at one end of the bursting regime. We therefore introduce a perturbative term in the slow equation, which extends this canonical form to a new one that we call \textit{Leidenator} and which supports the right transitions to and from elliptic bursting via classical and mixed-type torus canards, respectively. Throughout the study, we use singular flows ($\eps=0$) to predict the full system's dynamics ($\eps>0$ small enough). We consider three singular flows: slow, fast and average slow, so as to appropriately construct singular orbits corresponding to all relevant dynamics pertaining to elliptic bursting and torus canards.
\end{abstract}



\section{Introduction}
\label{sec:Introduction}
Complex oscillations are ubiquitous in neuronal data, where cells respond to sufficiently strong input stimuli not only by emitting action potentials or spikes but also by displaying more complicated electrical oscillations such as \textit{bursting}, which corresponds to an alternation between epochs of quiescent (subthreshold) activity and groups of consecutive spikes referred to as bursts. 
A large repertoire of bursting patterns are observed in neuronal experimental recordings~\cite{alving1968spontaneous,adams1985generation,butera1995analysis,cook1984electrical,jian2004novel,pinault1992voltage} as well as in other excitable cells~\cite{atwater1980nature}. Mathematical models of such complex neuronal oscillatory behaviour often possess multiple timescales, which can be explicit or not, and many different types of experimentally-observed bursting oscillations have been classified over the past few decades~\cite{rinzel1987formal,izhikevich2000neural,golubitsky01} using tools from bifurcation theory and singularity theory; see below. 

Biologists and modellers are interested in understanding possible changes (upon input) between different neuronal activity regimes, in particular how neurons may transition from spiking to bursting activity. One scenario of such transition involves special solutions called \textit{torus canards}; they have been evidenced and studied in both biophysical~\cite{kramer2008new} and phenomenological~\cite{benes2011elementary} models. According to this scenario, the membrane potential response to stimuli evolves from repetitive spiking (regular large-amplitude oscillations) to amplitude-modulated spiking, where a second frequency emerges in the solution profile, to then bursting activity. Torus canards are organising this transition, which occurs in a very narrow range of input parameter values.

Mapping the transition between bursting and tonic firing, which may carry
different information content~\cite{latorre2006,welday2011}, is a key question in neuroscience and it has
partially been addressed in single cell recordings\cite{ranck1973studies,
mason1988passive, connors1990intrinsic}, and investigated in a
FitzHugh–Nagumo–Rinzel type model \cite{wojcik2011}. As hinted at above, a first classification of bursting behaviour was proposed by Rinzel~\cite{rinzel1987formal} and later extended by Izhikevich
\cite{izhikevich2000neural}. These classifications hinge upon the slow-fast structure
of the models supporting bursting, and upon the bifurcation structure of certain
associated singular systems, which we now introduce.

In this paper, we will study \textit{elliptic bursting}~\cite{azad2010,ermentrout2005,lajoie2011,su2003effects}, which occurs
minimally in systems of $2$-fast $1$-slow variables. These models are written in
the \emph{fast-time parametrisation}
\begin{align}
  x' & =  f (x,y,\mu,\eps), && y' = g(x,y,\mu,\eps), && \mu' = \eps h(x,y,\mu,\eps), && t \in \mathbb{R}^+,
\label{eq:fastTime} \\
\intertext{or in an equivalent slow-time parametrisation}
  \eps \dot x & =  f (x,y,\mu,\eps), && \eps \dot y = g(x,y,\mu,\eps), && \dot \mu = h(x,y,\mu,\eps), && \tau \in \mathbb{R}^+,
\label{eq:slowTime}
\end{align}
where $x,y$ are fast variables, $\mu$ is the slow variable, $0 < \eps \ll 1$ is a
small parameter indicating the time scale separation, and in which prime and overdot
denote differentiations with respect to the fast time $t$ and to the slow time $\tau
= \eps t$, respectively. The functions $f,g,h$ are supposed to be sufficiently
regular.

As in many bursting patterns, the timescale difference in elliptic bursting systems
manifests itself as a composition of fast oscillatory components coupled to a slow
process driving the fast variables alternatively from quasi-stationary to
quasi-periodic dynamics~\cite{rinzel1987formal}. Intuitively, this means that one
uses both scalings presented above to understand bursting solutions.
Despite being equivalent to each other for
$\eps\neq 0$, systems~\cref{eq:fastTime} and~\cref{eq:slowTime} converge to two
non-equivalent singular limits for $\eps=0$:
\begin{align}
  x' & = f(x,y,\mu,0),  & & y' = g(x,y,\mu,0) & & \mu' = 0,                 & & t \in \mathbb{R}^+,  \label{eq:fastSubs}\\
   0 & = f(x,y,\mu,0),  & & 0  = g(x,y,\mu,0) & & \dot \mu = h(x,y,\mu,0),  & & \tau \in \mathbb{R}^+, \label{eq:slowSubs}
\end{align}
which are referred to as \textit{fast} and \textit{slow
subsystems}, respectively.

An important object in slow-fast dynamical systems is the \textit{critical manifold}
$S_0$
\[
  S_0:=\big\{(x,y,\mu)\in \mathbb{R}^3 \colon 0=f(x,y,\mu,0), \quad 0=g(x,y,\mu,0) \big\}.
\]
Note that $S_0$, which is a $1$-dimensional manifold in this case, appears both as
the set of equilibria of the fast subsystem (where the slow variable $\mu$ is frozen
and considered as a bifurcation parameter) and as the phase space of the slow
subsystem.

According to Izhikevich's classification, elliptic bursting is characterized by the
following pair of fast-subsystem bifurcations: a subcritical Hopf bifurcation (HB),
which triggers the onset of burst, and a saddle node (SN) of limit cycles, which
terminates the burst. This type of bursting is one example where the transition to and from bursting is associated with torus canards (TCs), which emerge in a exponentially-small parameter intervals~\cite{kramer2008new,benes2011elementary,burke2012showcase,ju2018bottom,vo2017generic,roberts2015averaging}.

We review these concepts in a neuronal population Wilson--Cowan-type model~\cite{wilson1972excitatory} which supports
elliptic bursting behaviour; see~\cite{burke2012showcase,izhikevich2000neural}. In the fast scaling, the model reads
\begin{equation}\label{eq:slowSS_meanField}
\begin{aligned}
x'= & -x+N_x(J_{xx}x+J_{xy}y+\mu) ,\\
y'= & \delta(-y+N_y(J_{yx}x+J_{yy}y+\rho)),\\
\mu^\prime= &\eps(k-x-y),
\end{aligned}
\end{equation}
where $x$ and $y$ model a fast excitatory and inhibitory population, respectively,
and where $\mu$ denotes a slow, excitatory, adaptation current. In the model
above, $N_x,N_y$ are the sigmoidal functions, 
\begin{equation}
N_x = \frac{\lambda_x}{2}\left[+\mathrm{erf}\left(\frac{g_1x}{\sqrt{2(1+g_1^2\sigma_x^2)}}\right)\right],\quad 
N_y = \frac{\lambda_y}{2}\left[1+\mathrm{erf}\left(\frac{g_2y}{\sqrt{2(1+g_2^2\sigma_y^2)}}\right)\right],
\end{equation}
and $J_{xx}$, $J_{xy}$, $J_{yx}$, $J_{yy}$, $\delta$, $\eps$, 
$\rho$, $g_1$, $g_2$, $v$, $\lambda_x$, $\lambda_y$, $k$, are control parameters; see also~\cite{touboul2015noise}.

\begin{figure}
\centerline{\includegraphics[width=0.8\textwidth]{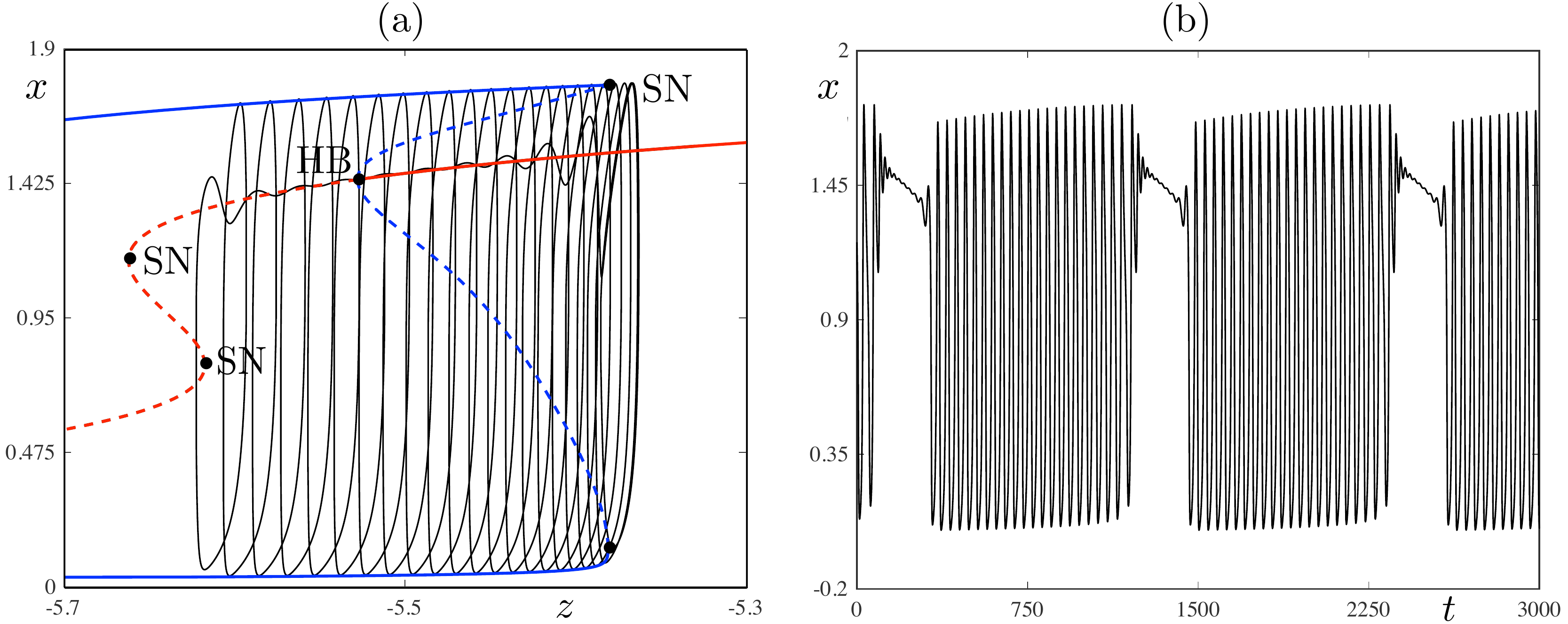}}
\caption{Elliptic bursting cycle from the Wilson-Cowan-type system~\cref{eq:slowSS_meanField}. (a) Projection onto the $(z,x)$ plane superimposed onto the bifurcation diagram of the fast subsystem. (b) Time profile of the $x$ variable over a few periods of the bursting cycle. Parameter values are: $J_{xx}=12$, $J_{xy}=-4$, $J_{yx}=13$, $J_{yy}=-9$, $\delta=0.05$, $\eps=0.001$, $\rho=-1.1$, $g_1=0.19$, $g_2=0.4$, $\sigma_x=\sigma_y=1.001$, $\lambda_x=2$, $\lambda_y=8.5$, $k=2.56$.}
\label{fig:ellipticClassicalScn}
\end{figure}

In~\cref{fig:ellipticClassicalScn} (a), we plot a bifurcation diagram of the fast subsystem, onto
which we superimpose an elliptic bursting cycle, whose time profile is shown in the panel (b) of the figure. The critical manifold $S_0$ is a cubic-like curve, and the
fast subsystem bifurcations relevant to elliptic bursting are the right-most HB and
the SN of limit cycles.

\begin{figure}
\centering
\includegraphics[width=\textwidth]{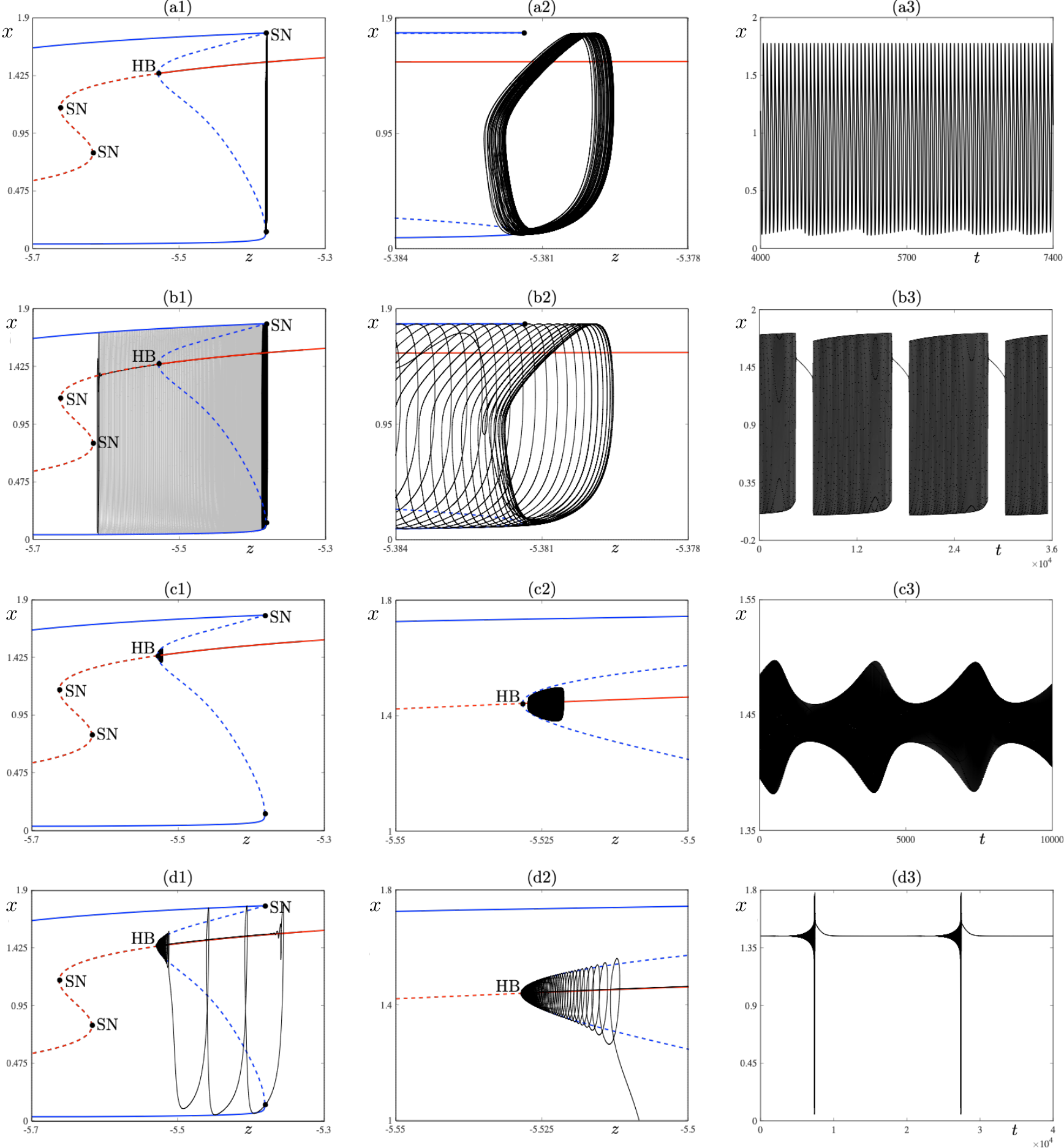}
\caption{Torus canard dynamics in the Wilson-Cowan type model~\cref{eq:slowSS_meanField}. Left and central panels show torus canard orbits in the $(z,x)$ plane superimposed onto the fast subsystem bifurcation diagram; the latter are zoomed views of the former. Right panels show the associated time series for $x$. Rows (ai), (bi), (ci) and (di) correspond in the same order to a headless TC, a TC with head, a headless mixed-type TC and a mixed-type TC with head. Fixed parameter values are as in~\cref{fig:ellipticClassicalScn}. Panels (ai)-(bi): $\eps=0.0001$, $k=2.536083$ (a) and $k=2.536084$ (b). Panels (ci)-(di): $\eps=0.001$, $k=3.60118$ (c) and $k=3.6011$ (d).}
\label{fig:TorCanWC}
\end{figure}

As parameter $k$ varies, we observe transitions to different regimes. When
$k$ decreases, we observe a transition to tonic
spiking. This transition is organised by TC cycles (see~\cref{fig:TorCanWC}(a)-(b)).
One geometrical property of TCs is that their fast oscillations do not stop
near the SN point of the fast subsystem, but instead continue along the branch of
unstable limit cycles. This feature justifies the term \textit{canard} because the
fast oscillations follow a repelling object of the fast
subsystem~\cite{BCDD1981,DR96,krupa01}. The name TC reflects the fact that this
transition is accompanied by a torus bifurcation in the full system (not shown), that
is, when $\eps\neq 0$ and $\eps\ll 1$ . As other canard solutions, torus canards
come both with head (\cref{fig:TorCanWC}(b)) and without head (\cref{fig:TorCanWC}(a)): the latter winds around an invariant torus of
the full system ($0 < \eps\ll 1$). 

In addition to the aforementioned one, the Wilson--Cowan-type system generates
another type of torus canard that seems to be less studied. We refer to it as \textit{mixed-type torus canard} for its link to the notion of canards of mixed type reported
in~\cite{desroches2012canards}. Like the classical TCs, these cycles are related to
the presence of a Neimark Sacker bifurcation in the full system. The name mixed type is used
because these solutions contain an attracting slow segment followed by a repelling
fast segment, hence they differ from classical TCs; mixed-type TCs appear near a
subcritical HB as shown in~\cref{fig:TorCanWC}(c)-(d). This is different from the classical TCs, which appear near a SN of cycles of the
fast subsystem.

In the present paper, we are motivated by studying minimal elliptic bursters that can
display all the solution types mentioned above. Classical TCs have been studied in
various neuronal models~\cite{kramer2008new,benes2011elementary,burke2012showcase}.
Mixed-type TCs, have received less attention, albeit they have been observed and
described in a recent study by Ju, Neiman and Shilnikov~\cite{ju2018bottom}, where
the main focus is on the full-system torus bifurcation.

The canonical form proposed by Izhikevich~\cite{izhikevich2000subcritical} for Bautin bursters, was adapted in~\cite{ju2018bottom} to the framework of elliptic bursting systems. Dynamics relevant to elliptic bursting, in particular torus canards, was therein studied with this adapted canonical form. We will refer to this adapted form as the \textit{canonical elliptic burster}. This model is considered to be a minimal
system for general elliptic bursting dynamics but, as we will show below, it does not
support mixed-type TCs. We first study this model and in particular we find that one boundary of the elliptic bursting regime is degenerate, owing to the presence of a continuum family of
equilibria. At such boundary, the critical manifold $S_0$ has a non-transversal
intersection with the slow nullcline. To overcome this limitation, we introduce a minimal perturbation of the canonical elliptic burster, which restores transversality at the aforementioned boundary. Consequently, the boundary in this new model is organised around mixed-type TCs. And away from this boundary, the new model behaves just like the canonical elliptic burster. Therefore, this new model appears to be a minimal torus canard system capturing transitions to and from the elliptic bursting regime.

In order to analyse all possible dynamics of the canonical elliptic burster as well as our minimal model, we mostly use singular flows and their concatenation, which define singular cycles. Full system's solutions, both of TC and of elliptic bursting type, stay close to such singular cycles for $\eps \neq 0$ small enough. While studying the different singular limits to construct such cycles has been a classical tool in planar canard systems~\cite{krupa01} as well as in systems with mixed-mode oscillations~\cite{brons2006}, it has not been used in bursting systems so far. Indeed, singular cycles perturbing to bursting cycles have not been studied at all, let alone in the TC regime; this is one of the main contributions of the present work.
To do so, we use an additional subsystem obtained in the
singular limit $\eps=0$, the so-called \textit{average slow subsystem}. This subsystem
approximates the slow flow on the manifold of limit cycles of the fast subsystem by averaging out the fast oscillations.  As a result,
the average slow flow retains a slowly-varying observable from the drifting fast
oscillations of the burst, and it allows to trace segments of singular solutions in
an additional $\eps=0$ slow limit. We then can realise all possible singular cycles by suitably concatenating segments from fast, slow and average slow subsystems, and we can understand how to recover singular equivalents of both classical and mixed-type TCs.

The manuscript is organized as follows. In \cref{sec:singular}, we introduce the concept of singular orbits and explain how they help to understand all possible dynamics of the full system, using the classical van der Pol model. Then in \cref{sec:canonicalforms}, we present and analyse both the canonical elliptic burster and our proposed generalized form; in particular, we focus on their singular orbits and relate them to full system's dynamics. Finally, in the conclusion section, we summarise our findings and propose a number of perspectives and questions for future work.

\section{Singular orbits}\label{sec:singular}
We review here the concept of singular
orbit~\cite{benoit1990,brons2006,krupa01}, which will be relevant for the analysis of
elliptic bursters presented below. The main idea behind singular orbits is to
extract information on solutions to the slow-fast dynamical system~\cref{eq:fastTime}, with
$0 < \eps \ll 1$, by constructing orbits for $\eps=0$. The latter are therefore
candidate orbits, that are hopefully descriptive of the perturbed orbits with $\eps
\neq 0$. As we have seen, there exist two non-equivalent singular 
limits in an $m$-slow $n$-fast dynamical system, namely the
fast subsystem~\cref{eq:fastSubs} and the slow subsystem~\cref{eq:slowSubs}, respectively. It is therefore natural to construct a singular orbit by concatenating
orbit segments of the fast- and slow-subsystem, as the following definition suggests
\begin{definition}[Singular orbit] \label{def:singOrbit}.
  Let $0=t_0 < t_1 < \cdots < t_l = T$
  be a partition of $[0,T] \subset \mathbb{R}$, for some $T>0$. A singular orbit
  starting at $p$ is a continuous curve $\gamma \colon [0,T] \to \mathbb{R}^{n+m}$
  \[
    \gamma(t) =
    \begin{cases}
      \gamma_1(t)  & t \in [t_0,t_1], \\
      \gamma_2(t)  & t \in [t_1,t_2], \\
      \ldots       &                 \\
      \gamma_{m}(t)  & t \in [t_{m-1},t_m]
    \end{cases}
    \qquad
    \gamma_{1}(t) = \varphi_1^t(p),
    \qquad
    \gamma_{i}(t) = \varphi_i^t(\gamma_{i-1}(t_{i-1})),
    \qquad
    i = 2,\ldots,l,
  \]
  where, for each $i=1,\ldots,l$, $(t,p) \mapsto \varphi_i^t(p)$ is the flow associated to either
  \cref{eq:fastSubs} or \cref{eq:slowSubs}.
\end{definition}

\begin{figure}
  \includegraphics[width=\textwidth]{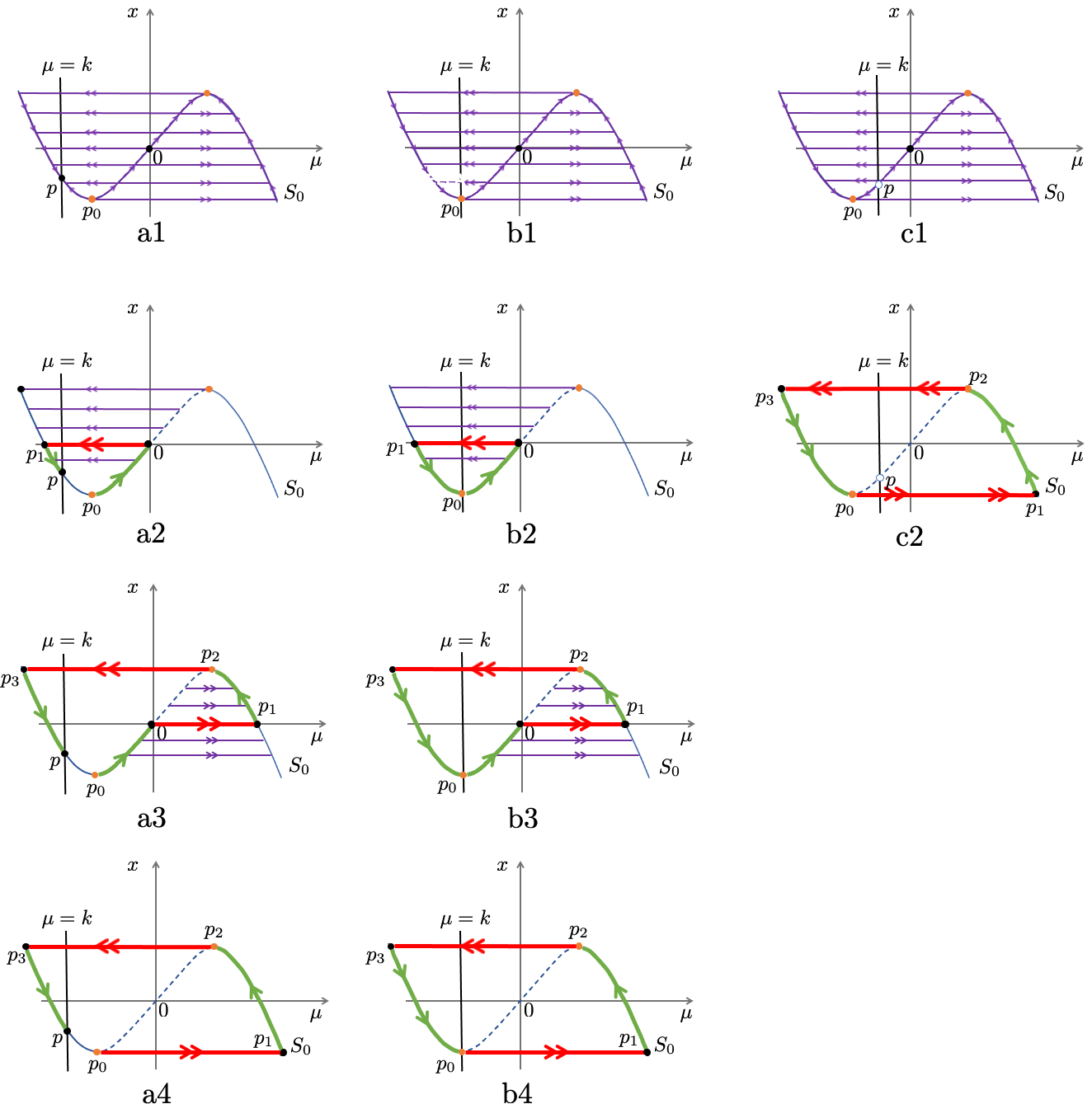}
  \caption{Examples of singular orbits in the Van der Pol system~\cref{eq:vandp} for
    different values of $k$ for which the slow nullcline intersects the critical
    manifold to the left of $p_0$ (a1--a4), at $p_0$ (b1--b4), and to the
    right of $p_0$ (c1--c2). (a1), (b1), (c1): sketches of slow (single arrow) and
    fast (double arrows) flows associated to~\cref{eq:vandp}, which give indication
    of how one can construct singular orbits. (a2)--(a4): Singular orbits starting at
    $p_0$ can follow the repelling branch of the critical manifold $S_0$ and then
    jump to the left (a2) or to the right (a3) attracting branch; we indicate slow
    segments in green, and fast segments in red; singular orbits are non unique, and
    we sketch in purple alternative fast segments; in (a4) we show an orbit starting
    at $p_0$ and jumping directly to the right attracting branch $S_0$; this singular
    orbit must follow $S_0$ up to the fold at $p_2$, from where it jumps towards
    $p_3$: it is impossible to concatenate to a fast segment before reaching the fold
    at $p_2$, while respecting \cref{def:singOrbit}, and indeed no purple segment is
    indicated. (b2)--(b4): Singular cycles can be constructed when the slow nullcline
    intersect $S_0$ at $p_0$, and similar considerations to (a2)--(a4) apply here;
    when $\eps >0$, these singular orbits perturb to headless canard cycles (a2),
    canard cycles with head (a3), and relaxation oscillations, respectively. (c2):
    when the $w$-nullcline intersects $S_0$ to the right of $p_0$ one can construct a
    unique singular cycle corresponding to a relaxation oscillation.
}
  \label{fig:vandP}
\end{figure}

In the analysis below, we will also be interested in periodic solutions to the full
problem, and hence, in closed singular orbits. We exemplify the construction of open
and closed singular orbits using the classical Van der Pol system~\cite{BCDD1981,DR96,krupa01}, known to be a prototypical planar canard system
\begin{equation}\label{eq:vandp}
  x' = x - \frac{x^3}{3} - \mu, \qquad \mu' = \eps(x-\mu).
\end{equation}
In \cref{fig:vandP} we show examples of singular trajectories. The top row sketches
slow and fast flows associated to~\cref{eq:vandp}, using single and double arrows,
respectively. The sketches are made for different values of $k$, for which the slow nullcline
intersects the critical manifold to the left of $p_0$ (a1), at $p_0$ (b1), and to the
right of $p_0$ (c1).

When the intersection is to the right of $p_0$, an equilibrium $p$ of the full
problem exists to the left stable branch of $S_0$. From panels (a2)-(a4) we see that
\emph{singular orbits are not unique}: concatenating slow (green) and fast (red) segments
and respecting the orientation on curve patches, one can construct infinitely many
singular orbits starting at
$p_0$, and terminating at $p$. In (a2) we see that an initial slow segment starting
from $p_0$ can be concatenated to any of the purple
segments, giving rise to orbits that jump to the left attracting branch of $S_0$, and
terminate at $p$. Of all such orbits, we highlight one that leaves the repelling
branch of $S_0$ at the origin (red segment). In (a3) we present a singular orbit that
starts at $p_0$ but then jumps at the origin to the right attracting branch of $S_0$,
highlighting in purple other potential concatenations. In (a4) we present a singular
orbit that starts at $p_0$, from where it jumps to the right attracting branch of
$S_0$ at $p_1$; in this case, the singular orbit must follow $S_0$ up to the fold at
$p_2$, from where it jumps towards $p_3$: it is impossible to concatenate to a fast
segment before reaching the fold at $p_2$, while respecting \cref{def:singOrbit}.

When the intersection of the $w$-nullcline is exactly at $p_0$, one can construct
singular cycles (panels (b1)--(b4)). Singular cycles may start at $p_0$, follow the
repelling branch and then jump to the left (b2) or to the right (b3) attracting
branch of $S_0$. A singular cycle that jumps from the starting point $p_0$ to the
right attracting branch of $S_0$ can also be uniquely constructed. When $\eps >0$,
the singular solutions give rise to headless canard cycles (b2), canard cycles with
head (b3), and relaxation a unique relaxation oscillation (b4).

Finally, when the intersection of the $w$-nullcline is to the right of $p_0$, one can
also construct a unique singular cycle corresponding to a relaxation oscillation.

\FloatBarrier

\section{Canonical forms for elliptic bursting}\label{sec:canonicalforms}
We now proceed to study the two canonical systems for elliptic bursting of interest
to us. The first one is the canonical elliptic burster in~\cite{ju2018bottom}, adapted  from~\cite[Thm 2.1]{izhikevich2000subcritical},
which we repeat here for convenience
\begin{equation}\label{eq:normalFormClassicalRep}
 \begin{aligned}
   x'   &=  -y+x\big ( \mu+2(x^2+y^2)-(x^2+y^2)^2  \big ),\\
   y'   &=   x+y\big ( \mu+2(x^2+y^2)-(x^2+y^2)^2 \big ),\\
 \mu' & = \eps (k-x^2-y^2).
 \end{aligned}
 \end{equation}
 A polar coordinate transformation $(x,y) = (r \cos \theta, r \sin \theta)$ leads to
 the system
\begin{equation}\label{eq:canonicalStandard}
\begin{aligned}
r' &=  r(\mu+2r^2-r^4),\\
\theta' &= 1,\\
\mu'  &=  \eps (k-r^2).
\end{aligned}
\end{equation}
We will consider $3$ subsystems for the canonical elliptic burster: in addition to the fast and
slow subsystems of \cref{eq:normalFormClassicalRep}, written as
\begin{equation}\label{eq:normalFormClassicalRepFast}
 \begin{aligned}
   x'   &=  -y+x\big ( \mu+2(x^2+y^2)-(x^2+y^2)^2  \big ),\\
   y'   &=   x+y\big ( \mu+2(x^2+y^2)-(x^2+y^2)^2 \big ),\\
 \mu' & = 0,
 \end{aligned}
 \end{equation}
 and
\begin{equation}\label{eq:normalFormClassicalRepSlow}
 \begin{aligned}
   0   &=  -y+x\big ( \mu+2(x^2+y^2)-(x^2+y^2)^2  \big ),\\
   0   &=   x+y\big ( \mu+2(x^2+y^2)-(x^2+y^2)^2 \big ),\\
   \dot \mu  &=  (k-r^2),
 \end{aligned}
 \end{equation}
 respectively, we will consider the averaged slow subsystem~\cite{butera1997transient,roberts2015averaging}, which is
 obtained by averaging the slow equation $\mu$ in the slow-time parametrisation of
 \cref{eq:normalFormClassicalRep} over one period of the fast limit
 cycles of the fast subsystem~\cref{eq:normalFormClassicalRepFast}, and approximating
 $\mu$ by its average $\langle \mu \rangle$. The fast-subsystem limit
 cycles, which we denote $(x_p(t;\mu),y_p(t;\mu))$, have period $2\pi$ for every
 $\mu>0$, hence we obtain
 \begin{equation}
   \dot{\langle \mu \rangle} = \frac{1}{2\pi} \int_0^{2\pi} \big( k - x^2_p(t;\langle
   \mu \rangle) - y^2_p(t;\langle \mu \rangle) \big)\, dt.
 \end{equation}

The fast, slow, and averaged slow subsystems discussed above have a compact
representation in polar coordinates, given by
\begin{align}
  r' &=  r(\mu+2r^2-r^4), && \theta' = 1, && \mu'  =  0 , \label{eq:BautinFS}\\
  0  &=  r,               &&              && \dot \mu  = k-r^2, \label{eq:BautinSS}\\
  0  &=  \mu +2r^2 - r^4, &&              && \dot \mu  = k-r^2, \label{eq:BautinASS}
\end{align}
respectively, where we have replaced $\langle \mu \rangle$ by $\mu$ in the last equation.
Note that there is no $\theta$ dynamics in none of the slow and average slow subsystems.

As stated in the introduction, in addition to the canonical elliptic burster, we
consider an extended model, with a perturbative term in the slow equation for $\mu$
\begin{equation}\label{eq:leidenator}
 \begin{aligned}
   x'   &=  -y+x\big ( \mu+2(x^2+y^2)-(x^2+y^2)^2  \big ),\\
   y'   &=   x+y\big ( \mu+2(x^2+y^2)-(x^2+y^2)^2 \big ),\\
   \mu'  &=  (k-x^2-y^2 - \alpha \mu),
 \end{aligned}
 \end{equation}
 where $\alpha >0$. We will henceforth refer to the model above, or to its polar
 representation
\begin{equation}\label{eq:LeidenStandard}
\begin{aligned}
r' &=  r(\mu+2r^2-r^4),\\
\theta' &= 1,\\
\mu'  &=  \eps (k-r^2 - \alpha \mu),
\end{aligned}
\end{equation}
 as the Leidenator model. A stochastic version of this model was simulated in a study on
 epilepsy~\cite[Supplementary material]{hebbink2017}. Reasoning as for the canonical
 elliptic burster, we arrive at the following slow, fast, and averaged slow
 subsystems for the Leidenator
\begin{align}
  r' &=  r(\mu+2r^2-r^4), && \theta' = 1, && \mu'  =  0 , \\
  0  &=  r,               &&              && \dot \mu  = k-r^2 - \alpha \mu, \\
  0  &=  \mu +2r^2 - r^4, &&              && \dot \mu  = k-r^2 - \alpha \mu,
\end{align}
respectively. As expected, the canonical and the Leidenator bursters differ only in
the slow subsystems.

In the following, we will analyse for each $k$ interval both the standard canonical
form and the Leidenator model in terms of singular and perturbed (non-singular)
dynamics. We will focus on selected singular orbits in each case, which are the most
relevant for the bursting.

\begin{figure}[!b]
\centerline{\includegraphics[scale=0.29]{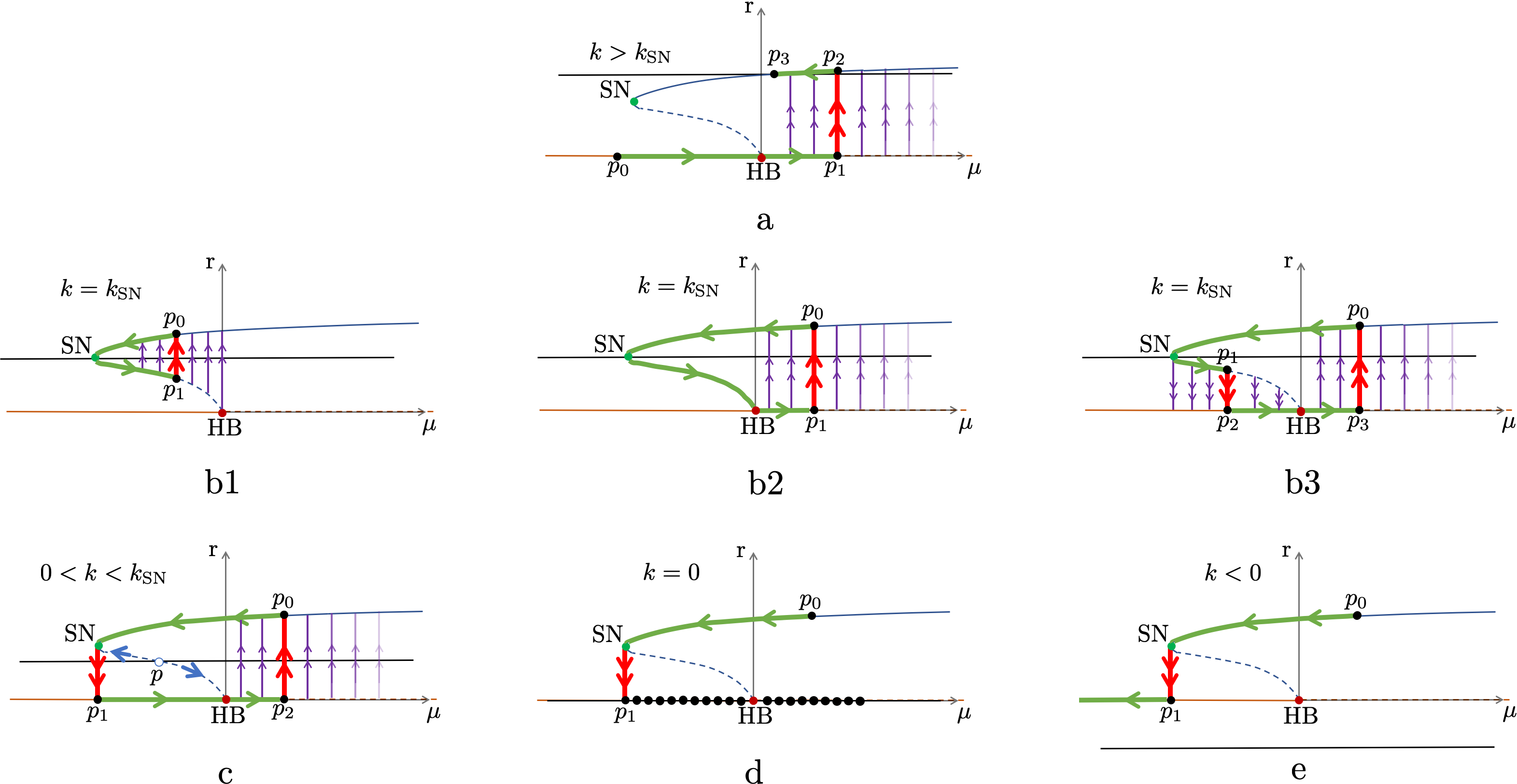}}
\caption{Sketches of the singular orbits generated by the canonical burster as $k$
  varies. Orbit segments of the slow/average slow subsystems
  \cref{eq:BautinSS,eq:BautinASS} are indicated in green, and segments of the fast
  subsystem \cref{eq:BautinFS} are indicated in red. We denote by SN and HB
  bifurcations of the fast subsystem. Singular orbits are constructed by
  concatenating red and green segments, with concatenation points indicated by dots.
  Purple fast segments indicate alternative fast orbit segments. Solid (dashed)
  curves indicate branches of stable (unstable) attractors of the fast subsystem.
}
\label{fig:preLeidenator}
\end{figure}

\subsection{Canonical elliptic burster}
We now consider the fast subsystem \cref{eq:BautinFS} of the canonical elliptic
burster. The critical manifold is given by the horizontal line $r=0$. The fast
subsystem has 2 main bifurcations, a subcritical Hopf, HB, and a saddle-node of
cycles, SN, as shown in \cref{fig:preLeidenator}. The fast subsystem is therefore
bistable, and the branch of periodic orbits emanating from HB are circles of radius
$r = \sqrt{\mu}$, and lie in the $(r,\mu)$ plane on the curve $0 = \mu + 2r^2 -r^4$.

We will study the dynamics of the full system when $k$ is varied. As we shall see
below, the asymptotic dynamics of the full system is either a stable limit cycle, or
a bursting cycle, and the transition between these regimes occurs around a value $k =
\sqrt{r_{SN}}$, where $r_{SN}$ is the fold point of $0=\mu + 2r^2 - r^4$.
Henceforth, we will pose $k_{SN} =\sqrt{r_{SN}}$. We will also show below that the bursting regime
terminates around $k=0$, where the dynamics is degenerate.

\subsubsection*{Case 1: $\boldsymbol{k>k_{SN}}$}
\paragraph{Singular dynamics} The slow average subsystem has a unique stable
equilibrium at $p_3$ as shown in~\cref{fig:preLeidenator}(a). All singular orbits terminate at $p_3$ in this region of parameter space. We show
the construction of a typical orbit, which starts on the critical manifold at $p_0$,
follows the critical manifold past HB, up to an arbitrary large value $p_1$. From
$p_1$ the orbit concatenates with a fast segment to $p_2$; the final segment from
$p_2$ to $p_3$ is on the slow averaged flow. There is a one-parameter family of
orbits connecting $p_0$ to $p_3$; this family is parametrised by $p_1$.

\paragraph{Nonsingular dynamics}
The full system ($\eps \neq 0$) selects a unique orbit starting at $p_0$; this orbit
approaches an $\eps$-perturbation of $p_3$, that is, a tonic firing solution of the
full system. This orbit lifts off from an $\eps$-neighbourhood of the critical manifold, at a
point uniquely determined using the so-called \textit{way-in-way-out
function}~\cite{diener1991} once $k$, $\eps$ and $p_0$ are fixed.

\subsubsection*{Case 2: $\boldsymbol{k=k_{\text{SN}}}$}

\paragraph{Singular dynamics}
In this case, the $\mu$-nullcline crosses the fold at SN, hence the slow average system does not have an equilibrium, but a turning
point at SN and, as for the van der Pol system, this generates singular canard cycles,
as presented in the three panels \cref{fig:preLeidenator}(b1)--(b3).

In the first scenario, which is illustrated in~\cref{fig:preLeidenator}(b1), we
construct a one-parameter family of singular orbits parametrized by $p_1$,
consisting of one average slow segment, which goes through SN, and
one fast segment. SN is a canard point for the $(r,\mu)$ averaged slow flow
system~\cref{eq:BautinASS}. Each member of this one-parameter family is uniquely defined based
on the location of $p_1$, which can vary between SN and HB (excluded), thereby generating
singular canard orbits in~\cref{eq:BautinASS}. We call these orbits singular TC
orbits because they correspond to TC orbits in the full system.

In the second scenario, in~\cref{fig:preLeidenator}(b2), $p_1$ occurs at HB, or at an
arbitrary point to the right of HB on the critical manifold. These orbits are the
so-called \textit{singular maximal TCs}, because their repelling segment on
the branch of periodic orbits of the fast subsystem has maximal length.

In the third scenario we observe a
two-parameter family of singular TCs with head (with 2 fast segments) and parametrized by $p_1$ and $p_3$;
see \cref{fig:preLeidenator}(b3). The location of $p_1$, which is arbitrarily chosen
between SN and HB, determines the canard's head size. As in the second scenario,
$p_3$ can be located anywhere on the unstable critical manifold.

\paragraph{Nonsingular dynamics} When $\eps \neq 0$, we expect to find TC cycles
in the full system~\cref{eq:normalFormClassicalRep}; these orbits will be close to
singular orbits from the three scenarios mentioned above, and parametrised by $k$ in
an exponentially small region of parameter space. This one-parameter family of TC
cycles starts with headless TCs realised from~\cref{fig:preLeidenator}(b1). Out of
the singular maximum TC cycles of~\cref{fig:preLeidenator}(b2), we expect that only
one is realised, namely the one for which $p_1$ is at HB. This is best understood by
looking at the third scenario~\cref{fig:preLeidenator}(b3): computing the
way-in-way-out function, we see that $p_2$ and $p_3$ must be equidistant from HB.
This also implies that, when $\eps \neq 0$, the TC cycles with realised in the full
system are a one-parameter family, depending solely $p_2$, as opposed to their
singular counterpart which appear in a 2-parameter family parametrised by $p_2$ and
$p_3$.

\subsubsection*{Case 3: $\boldsymbol{0<k<k_{\text{SN}}}$}

\paragraph{Singular dynamics}
In this case the equilibrium of the averaged slow system~\cref{eq:BautinASS} is on
the unstable branch of limit cycles, therefore singular orbits do not display canard
segments, as they connect SN to a point $p_1$ on the critical manifold. One can
construct a one-parameter family of singular cycles, parametrised by the value of
$p_2$ at which the singular orbit is concatenated to a fast segment, as shown
in~\cref{fig:preLeidenator}(c). We call these orbits singular elliptic bursting
cycles.

\paragraph{Nonsingular dynamics}
For each $k,\eps>0$, we expect that only elliptic bursting orbit persists, selected by the way-in-out
function as above.

\subsubsection*{Case 4: $\boldsymbol{k=0}$}

\paragraph{Singular dynamics}
In this case, the slow subsystem has a continuum of trivial equilibria on the
critical manifold, as illustrated in~\cref{fig:preLeidenator}(d). All singular orbits
starting at a point $p_0$ on the stable slow average branch terminate at $p_1$, the
projection of SN on $S_0$. One can construct other singular orbits, terminating at
any of the equilibria on $S_0$, by varying $p_0$.

\paragraph{Non-singular dynamics}
$S_0$ is now also a line of equilibria for the full
system~\cref{eq:normalFormClassicalRepFast}. The non-singular behaviour is
qualitatively the same of the singular one: the system comes to a rest on $S_0$,
after a jump from the branch of stable limit cycles.

\subsubsection*{Case 5: $\boldsymbol{k<0}$}
\paragraph{Singular and nonsingular dynamics}
There is no equilibrium of the slow subsystems. Both the singular and the nonsingular
dynamics display a drift in $\mu$ towards $-\infty$, while $r(t)$ tends monotonically
to $0$.

\subsection{The Leidenator model}
The canonical bursting model captures well the transition from tonic firing to
elliptic bursting via TCs, observed in most known elliptic bursters, by decreasing
$k$. However, as $k$ decreases further, the end of the elliptic bursting regime
occurs via a continuum of equilibria, which is non generic, because the fast nullcine
(the critical manifold) intersects non transversally the $\mu$-nullcline (the two
coincide). In addition, in numerical simulations of the Wilson-Cowan-type model~\cref{eq:slowSS_meanField}, we have seen that the end of the bursting
regime occurs via mixed-type TC cycles, like the ones presented
in~\cref{fig:TorCanWC}(ci)-(di). The canonical elliptic burster can not capture
solutions of this type.

We introduce the Leidenator model in order to overcome these problems. In this model, the slow nullcline is the parabola $\mu = (k - r^2) / \alpha$ and it intersects transversally the critical manifold for all $k$ values. For sufficiently small $\alpha>0$, one therefore obtains an isolated equilibrium of the slow flow; see~\cref{fig:preLeidenator}. In addition, for $k>0$, this equilibrium corresponds to a so-called \textit{buffer point}~\cite{diener1991}, which affects the non-singular dynamics as will be explained below. Note that one recovers the canonical elliptic burster from the Leidenator for $\alpha=0$.

\subsubsection*{Cases 1-3: $\boldsymbol{k>0}$}
\begin{figure}
\centerline{\includegraphics[scale=0.27]{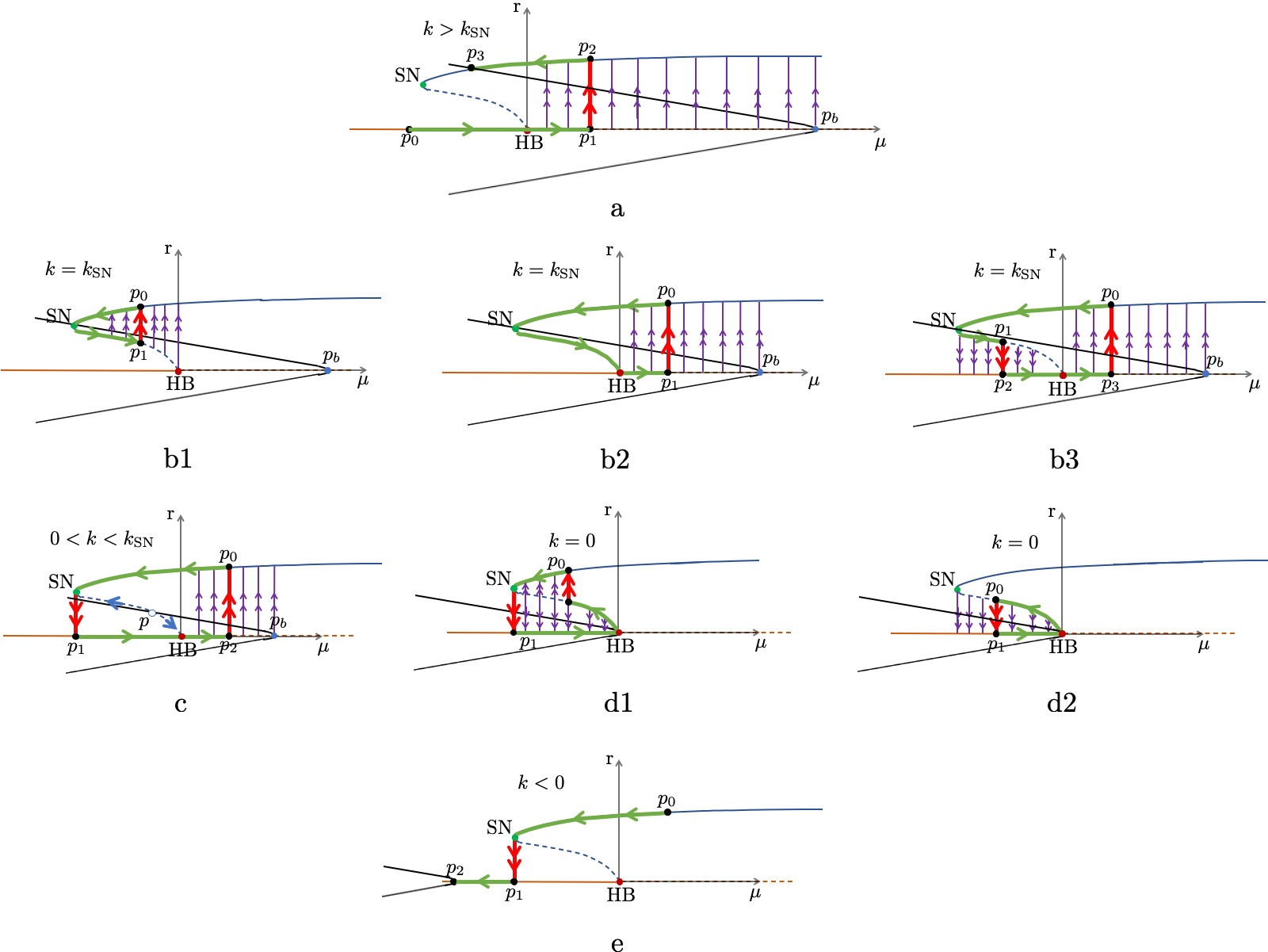}}
\caption{Sketches of the singular orbits generated by the Leidenator model. Graphical conventions as in \cref{fig:preLeidenator}. A buffer point $p_b$ arises and is indicated by a blue dot. Panels (d1) and (d2) show new configurations, which give rise to singular mixed-type TC with head and without head, respectively.}
\label{fig:captureLeidenator}
\end{figure}

As seen in \cref{fig:captureLeidenator}(a)-(c), the Leidenator behaves like the canonical elliptic bursting model for $k>0$. The only difference is the existence of $p_b$, which determines an upper bound for the points on $S_0$ where the right-most fast segment can be placed. For $\eps \neq 0$, the buffer point at $p_b$ marks the maximal delay that any full-system trajectory can have.
\subsubsection*{Case 4: $\boldsymbol{k=0}$}
\begin{figure}
\centerline{\includegraphics[scale=0.45]{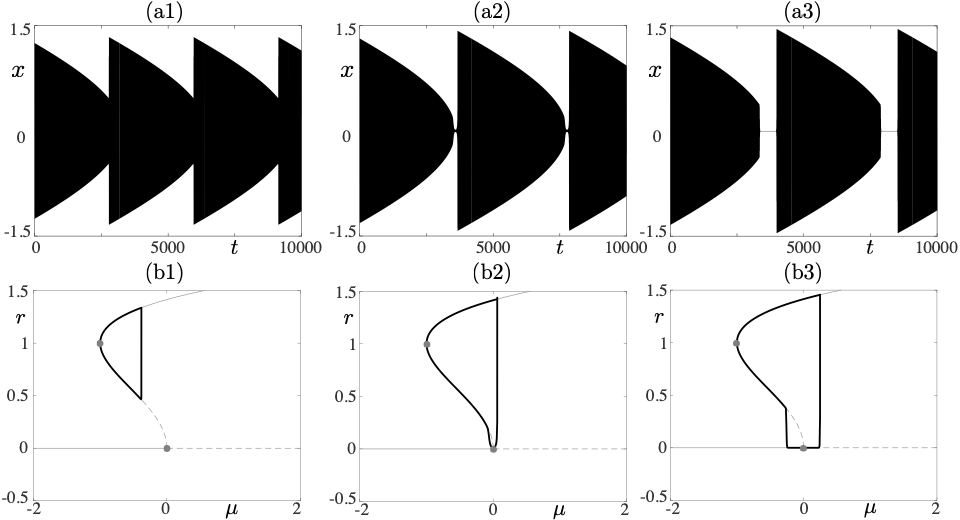}}
\caption{TC cycles computed in the Leidenator model. From left to right: headless TC, maximal TC and TC with head. Top panels display the time profile of the $x$-variable ($x=r\cos\theta$). Bottom panels show the projection of the solution onto the $(\mu,r)$ phase plane together with the fast subsystem bifurcation diagram (the dots correspond to the HB and the SN points). Solutions of these types are also supported by the canonical elliptic burster. Parameter values: $\eps=0.001$, $\alpha=0.2$. The value of $k$ is different in each panel with numerical values within $10^{-12}$ of  $0.79984990182$.}
\label{fig:torcanleiden}
\end{figure}
\paragraph{Singular dynamics}
This case differs substantially from the canonical elliptic burster. Indeed, $p_b$ coincides with HB and there is no equilibrium of the average slow flow (point $p$ from \cref{fig:preLeidenator}(c) is no longer present). This provides a new way to concatenate singular segments, namely slow and average slow segments, so as to form new types of singular canard cycles. We refer to these cycles as singular mixed-type TC cycles. Similar to all singular canard cycles, we come both with and without head, which are shown in~\cref{fig:captureLeidenator}(d1) and (d2), respectively. At the transition between these two subfamilies, there is a unique singular maximal mixed-type TC cycle.

\begin{figure}
\centerline{\includegraphics[scale=0.45]{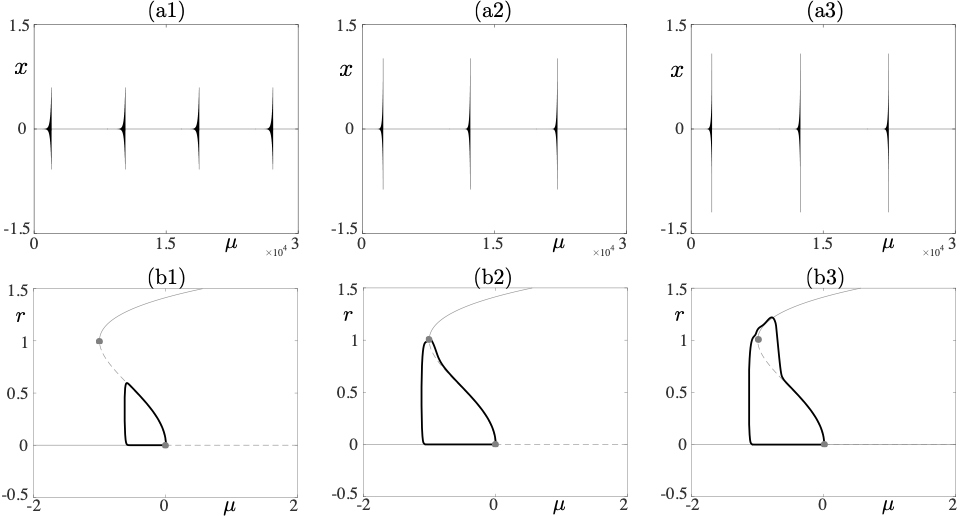}}
\caption{Mixed-type TC cycles computed in the Leidenator model. From left to right: headless mixed-type TC, maximal mixed-type TC and mixed-type TC with head. Top panels display the time profile of the $x$-variable ($x=r\cos\theta$). Bottom panels show the projection of the solution onto the $(\mu,r)$ phase plane together with the fast subsystem bifurcation diagram (the dots correspond to the HB and the SN points). These solutions are not supported by the canonical elliptic burster. Parameter values: $\eps=0.05$, $\alpha=0.2$. The value of $k$ is different in each panel with numerical values within $10^{-14}$ of  $0.0015128438002$.}
\label{fig:mixtyptorcanleiden}
\end{figure}

\paragraph{Non-singular dynamics}
When $\eps \neq 0$, we expect to find mixed-type TC cycles
in the full system~\cref{eq:leidenator}; these orbits will be parametrised by $k$ in
an exponentially small region of parameter space and will span canard cycles with and without head.

\subsubsection*{Case 5: $\boldsymbol{k<0}$}
This case is identical to Case 5 in the canonical elliptic burster.

\subsection{Numerical results}
We have used numerical continuation with the software package \textsc{auto}~\cite{auto} in order to compute canard cycles in the $(r,\mu)$ system and therefore obtain TC cycles in the Leidenator model via the polar coordinate transform. The results are included in \cref{fig:torcanleiden}, which shows TC cycles corresponding to Case 2, as well as in~\cref{fig:mixtyptorcanleiden}, which shows mixed-type TC cycles corresponding to Case 4.

\section{Conclusion}
\label{sec:conclusion}
In this paper, we analysed the dynamical repertoire of elliptic bursting systems by using two minimal systems, namely the canonical elliptic burster~\cref{eq:normalFormClassicalRep} and the Leidenator~\cref{eq:leidenator}. The presented results focuses on singular orbits in order to show that the canonical elliptic burster cannot produce all dynamical patterns that an elliptic burster may display. More precisely, it cannot possess mixed-type TCs due to the presence of a continuum of equilibria which correspond to one boundary of the elliptic bursting regime where mixed-type TC should appear.  For this reason, we extended the canonical elliptic burster into the Leidenator, which manage to overcome this limitation and possess generic transitions involving torus canards (classical and mixed type) at both ends of the bursting regime.

Our approach is based on the singular limits of slow-fast systems of elliptic bursting type, and it employs three limiting systems obtained for $\varepsilon = 0$: the fast, slow and average slow subsystems. We define a singular orbit as a compatible concatenation of orbit segments solutions to these three subsystems. The idea behind the singular orbits is to extract information about the full system with $0<\varepsilon \ll 1$ by constructing piecewise-defined orbits valid for $\varepsilon = 0$. The orbits of the full system remain close to the singular orbits. Therefore, the singular orbits contain information about the full system within a simplified setting.

In parameter space, the transition to and from elliptic bursting is organized by TCs, which appear in exponentially-small parameter intervals. Two families of TCs are then of interest: classical and mixed-type. In both families, canard solutions are divided into headless canards, maximal canard and canards with head. We relate such TC solutions, as well as the bursting and tonic spiking solutions, to appropriately constructed singular orbits. As we do so, the fact that the canonical elliptic burster~\cref{eq:normalFormClassicalRep} cannot produce mixed-type TCs becomes clear and show the limitation of this model, as mixed-type TCs are observed in elliptic bursters, e.g. in the Wilson-Cowan-type model given by~\cref{eq:slowSS_meanField}. This is due to the presence of a continuum of equilibria resulting from a nontransversal intersection between the (slow) $\mu$ nullcline and the critical manifold for $k=0$, as shown in~\cref{fig:preLeidenator}(d1). We then propose the Leidenator system, in which we introduce a perturbative term in the slow flow that minimally suffices to overcome this problem. This term makes the $\mu$ nullcline (a parabola) keep a transversal intersection with the critical manifold for any $k$ value, as shown in~\cref{fig:captureLeidenator}. Consequently, this guarantees that, for proper parameter choices, the full system has a unique equilibrium when $k=0$. At the level of singular systems, it enables the construction of singular mixed-type TCs when $k=0$, which is not possible in the case of the canonical elliptic burster.  

Using various $\eps=0$ systems in order to construct singular cycles, and hence infer the possible dynamics of a slow-fast system, is a standard tool that has been regularly used in problems with canard cycles~\cite{krupa01} or mixed-mode oscillations~\cite{brons2006}. However, it has hardly been used in the context of bursting systems, where a third singular flow must be considered. The present paper shows the advantage of this approach in order to derive valuable information in a given family of bursting systems, namely elliptic bursters. We anticipate that it could give useful in order type of bursters, starting by the other two classes introduced by Rinzel~\cite{rinzel1987formal} which are square-wave and parabolic bursters. This is an interesting question for future work. Also, the question of rigorous existence of mixed-type TCs is very natural and we will consider it in a follow-up study. Finally, 
the occurrence of mixed-type TCs in biophysical neuron model, and their physiological relevance, is also of interest, possibly in link with resonator neurons that display subthreshold oscillations, as well as the possibility to find such structures in a spatially-extended context.

\bibliographystyle{siamplain}
\bibliography{references}
\end{document}